\documentclass[a4paper,reqno]{amsart}

\usepackage{pdfsync}
\usepackage{amsmath,amssymb}
\usepackage{enumerate}
\usepackage{xspace}
\usepackage{boxedminipage}
\usepackage{algorithmic}
\usepackage{listings}% for including source code
\usepackage{tabularx}
\usepackage{subfigure}
\newtheoremstyle%
    {plain}% name
    {}% Space before, vuoto = `valore di default'
    {}% Space after
    {\mdseries\slshape}% body font
    {}% Indent (empty = no indent, \parindent = para indent)
    {\bfseries}% Thm head font
    {.}% Punctuation after the heading
    {1.0ex}% Space after heading: \newline = to start at next line
    {}% Thm head spec (can be left empty, meaning `normal')
    
\newtheoremstyle
    {note}% name
    {}% Space before, vuoto = `valore di default'
    {}% Space after
    {}% body font
    {}% Indent (empty = no indent, \parindent = para indent)
    {\bfseries}% Thm head font
    {.}% Punctuation after the heading
    {1.0ex}% Space after heading: \newline = to start at next line
    {}% Thm head spec (can be left empty, meaning `normal')
    
\swapnumbers
\theoremstyle{plain}
\newtheorem{The}[subsection]{Theorem}
\newtheorem{Pro}[subsection]{Proposition}
\newtheorem{Pro*}{Proposition}

\theoremstyle{note}

\newtheorem{Rem}[subsection]{Remark}
\newtheorem{Defn}[subsection]{Definition}

\usepackage{ifthen}
\usepackage{amsmath,amssymb}    %mathematics
\usepackage{stmaryrd} %more symbols
\usepackage{eufrak}     %fraktur charcters
\usepackage{mathrsfs} %nice calligraphic characters
\usepackage{bbold}    %blackboard bold
\provideboolean{usemathrsfs}
\setboolean{usemathrsfs}{true}
\usepackage{enumerate}% a facility to customize numbered lists
\usepackage{xspace}% smart spacing after macros: detects punctuation
\usepackage{boxedminipage}% mainly used for slides
\usepackage[latin1]{inputenc}% for accents, umlauts, etc. directly in the text
\usepackage{verbatim}% provides also the {comment} environment
\provideboolean{usebeamer}% when using beamertex, 
\provideboolean{isthesis}% used in the PhD thesis
\provideboolean{isamsltex}% use when using AMS papers style
\provideboolean{issiamltex}% on when using SIAM papers... becoming obsolete

\usepackage{hyperref}

\newcommand{\linkedemail}[1]{{\makeatletter\color{blue}\href{mailto:#1}{#1}\makeatother}}
\newcommand{\Ignore}[1]{}
\newcommand{\freeze}[1]{}% Obsolete, kept for compatibility
\newcommand{\crossout}[1]{{\textcolor{red}{#1}}}
\newcommand{\highlight}[1]{{\textcolor{blue}{#1}}}
\provideboolean{shownotes}
\setboolean{shownotes}{true}
\newcounter{margnote}[page]
\newcommand{\margnotemark}{\highlight{\upshape\texttt{>\arabic{margnote}<}}}
\newcommand{\margnote}[2][]{\ifthenelse{\boolean{shownotes}}%
{\stepcounter{margnote}%
\margnotemark%
\marginpar{\texttt{\raggedright\tiny\margnotemark{#1}: #2}}}
{}}
\provideboolean{showtodo}
\newcommand{\todo}[1]
{\ifthenelse{\boolean{showtodo}}{\margnote[To do.]{#1}}{}}
{\ifthenelse{\boolean{showtodo}}{\end{boxedminipage}}{}}

\provideboolean{showcomments}
\newcommand{\margincomment}[1]{
\ifthenelse{\boolean{showcomments}}{\marginpar{\tiny #1}}{}
}
\provideboolean{showchanges}
\setboolean{showchanges}{false}
\newcommand{\changes}[1]{
  \ifthenelse{\boolean{showchanges}}
	     {{\highlight{#1}}}
	     {#1}
}
\newcommand{\changefromto}[2]{
  \ifthenelse{\boolean{showchanges}}
  {{\crossout{#1}${\color{magenta}\mapsto}$}{\highlight{#2}}}
  {#2}
}

\newcommand{\ie}{\ensuremath{\text{ i.e., }}\xspace}

\usepackage{mathrsfs}%extra mathematical symbols
\provideboolean{usemathrsfs}
\setboolean{usemathrsfs}{true}% switch if mathrsfs is unused
\newcommand{\mathscript}
	   {\mathscr}
 \newcommand{\cE}{\ensuremath{\mathscript E}\xspace}
 \newcommand{\cF}{\ensuremath{\mathscript F}\xspace}

 \newcommand{\cT}{\ensuremath{\mathscript T}\xspace}

 \newcommand{\rN}{\ensuremath{\mathbb N}\xspace}
 
 \newcommand{\rP}{\ensuremath{\mathbb P}\xspace}
 
 \newcommand{\rR}{\ensuremath{\mathbb R}\xspace}
 \newcommand{\rS}{\ensuremath{\mathbb S}\xspace}

 \newcommand{\rW}{\ensuremath{\mathbb W}\xspace}

 \newcommand{\naturals}{\rN\xspace}

 \newcommand{\reals}{\rR}

 \newcommand{\closure}[1]{\overline{#1}}

 \newcommand{\W}{\ensuremath{\varOmega}\xspace}
 
 \newcommand{\qp}[1]{\ensuremath{\left({#1}\right)}}

 \newcommand{\powqp}[2]{\ensuremath{\qp{#2}^{\kern -.2em\lower .7ex\hbox{\scriptsize $#1$}}\kern-.3em}}
 
 \newcommand{\norm}[1]{\ensuremath{\left|#1\right|}}
 
 \newcommand{\Norm}[1]{\ensuremath{\left\|#1\right\|}}
 
 \newcommand{\ltwop}[2]{\ensuremath{\left\langle#1,#2\right\rangle}}

 \newcommand{\duality}[2]{\ensuremath{\left\langle #1\,\vert\,#2\right\rangle}}

 \renewcommand{\d}{\ensuremath{\,\mathrm{d}}}
 
 \providecommand{\D}{\ensuremath{\mathrm{D}}}

 \newcommand{\registered}%
	    {\ensuremath{{}^{\bigcirc\!\;\!\!\!\!\!\!\!\;\text{\sc r}}}}

\newcommand{\AND}{\ensuremath{\text{ and }}}

\renewcommand{\div}{\operatorname{div}}

\newcommand{\Oh} {\operatorname{O}}                   % Landau's Big O
\newcommand{\trace}{\operatorname{trace}}             % trace (of a matrix)
\newcommand{\transpose}{{\boldsymbol\intercal}}   % transpose symbol
\newcommand{\Transpose}[1]{\ensuremath{{#1}^{\transpose}}}

\newcommand{\Hess}{\ensuremath{\D^2}}
\newcommand{\intersected}{\ensuremath{\cap}}
\newcommand{\meet}{\intersected}

\newcommand{\union}[1]{\ensuremath{\bigcup}_{#1}}
\renewcommand{\vec}[1]{\ensuremath{\boldsymbol{#1}}}
\newcommand{\geovec}[1]{\ensuremath{\vec{#1}}}

 \newcommand{\CC}{\ensuremath{\operatorname C}\xspace}%Continuous functions
 \newcommand{\HH}{\ensuremath{\operatorname H}\xspace}
 \newcommand{\LL}{\ensuremath{\operatorname L}\xspace}

 \newcommand{\cont}[1]{\ensuremath{\CC^{#1}}}

 \newcommand{\leb}[1]{\ensuremath{\LL_{#1}}}

 \newcommand{\sobh}[1]{\ensuremath{\HH^{#1}}}
 \newcommand{\sobhz}[1]{\sobh{#1}_0}
 \newcommand{\poly}[1]{\ensuremath{\rP}^{#1}}

 \newcommand{\fes}[1]{\ensuremath{\fespace^{#1}}}
 \newcommand{\fesW}{\ensuremath{\rW}}
\newcommand{\Forall}{\:\forall\:}

\newcommand{\Foreach}{\quad\Forall}
\usepackage{amscd}

\newcommand{\funk}[3]{\ensuremath{#1:#2\to#3}}

\newcommand{\Program}[1]{\textsf{#1}\xspace}

\providecommand{\ListParameters}{}
\renewcommand{\ListParameters}% these can be changed by renewcommand
{
	 \setlength{\topsep}{0em}
	 \setlength{\leftmargin}{0em}
         \setlength{\itemsep}{0ex}
	 \setlength{\parsep}{.5ex}
	 \setlength{\itemindent}{\labelsep}
	 \addtolength{\itemindent}{\labelwidth}
}

\newcounter{LetterListItem}
\renewcommand{\theLetterListItem}{(\alph{LetterListItem})}
{
	\begin{list}%
	{\theLetterListItem\ }%
	{\usecounter{LetterListItem}
	 \ListParameters
	}
}%
{\end{list}}
\newcounter{NumberListItem}
\renewcommand{\theNumberListItem}{\arabic{NumberListItem}}
{
	\begin{list}%
	{\theNumberListItem.\ }%
	{\usecounter{NumberListItem}%
	 \ListParameters
	}
}%
{\end{list}}
\newcounter{QuestionListItem}
\renewcommand{\theQuestionListItem}{\textbf{Question \arabic{QuestionListItem}}}
{
	\begin{list}%
	{\theQuestionListItem.\ }%
	{\usecounter{QuestionListItem}%
	 \ListParameters
	}
}%
{\end{list}}
\newcounter{RomanListItem}
\renewcommand{\theRomanListItem}{(\roman{RomanListItem})}
{
	\begin{list}%
	{\theRomanListItem\ }%
	{\usecounter{RomanListItem}
	 \ListParameters
	}
}%
{\end{list}}
\newcounter{StepsItem}
{
	\begin{list}%
	{Step \theStepsItem.\ }%
	{\usecounter{StepsItem}%
	 \ListParameters
	}
}%
{\end{list}}

\usepackage{amsthm} % an extra package for theorem-like environments
\ifthenelse{\boolean{isthesis}}{%
  \setcounter{secnumdepth}{1}%
}{
  \setcounter{secnumdepth}{2} % this number is the depth of numbering
}
\providecommand{\ListParameters}{}
\renewcommand{\ListParameters}
{
	 \setlength{\topsep}{0em}
	 \setlength{\leftmargin}{0em}
         \setlength{\itemsep}{0ex}
	 \setlength{\parsep}{.5ex}
	 \setlength{\itemindent}{\labelsep}
	 \addtolength{\itemindent}{\labelwidth}
}
\newtheoremstyle{plain}% name
  {}% Space before, vuoto = `valore di default'
  {}% Space after
  {\mdseries\slshape}% body font
  {\parindent}% Indent (empty = no indent, \parindent = para indent)
  {\bfseries}% Thm head font
  {.}% Punctuation after the heading
  {.5em}% Space after heading: \newline = to start at next line
  {}% Thm head spec (can be left empty, meaning `normal')

\newtheoremstyle{note}% name
  {}% Space before
  {}% Space after
  {}% body font
  {\parindent}% Indent (empty = no indent, \parindent = para indent)
  {\bfseries}% Thm head font
  {.}% Punctuation after the heading
  {.5em}% Space after heading: \newline = to start at next line
  {}% Thm head spec (can be left empty, meaning `normal')

\newtheoremstyle{claim}% name
  {}% Space before, vuoto = `valore di default'
  {}% Space after
  {\mdseries\slshape}% body font
  {}% Indent (empty = no indent, \parindent = para indent)
  {\bfseries}% Thm head font
  {:}% Punctuation after the heading
  {.5em}% Space after heading: \newline = to start at next line
  {}% Thm head spec (can be left empty, meaning `normal')

\newtheoremstyle{exercise}% name
  {}%  Space before, empty = `default'
  {}% Space after
  {}% body font
  {}% Indent (empty = no indent, \parindent = para indent)
  {\bfseries}% Thm head font
  {.}% Punctuation after the heading
  {1em}% Space after heading: \newline = to start at next line
  {}% Thm head spec (can be left empty, meaning `normal')

\newtheoremstyle{break}% name
  {}%  Space before, empty = `default'
  {}% Space after
  {}% body font
  {}% Indent (empty = no indent, \parindent = para indent)
  {\bfseries}% Thm head font
  {.}% Punctuation after the heading
  {\newline}% Space after heading: \newline = to start at next line
  {}% Thm head spec (can be left empty, meaning `normal')

  \newcommand{\Proofname}{Proof}%Dimostrazione   %DÃÂÃÂÃÂÃÂ©monstration%Beweis

\ifthenelse{\boolean{isthesis}}{%
  
}{%
  
}

\newcommand{\pdfformat}[1]{
  \provideboolean{pdfoutput}
  \setboolean{pdfoutput}{#1}% default is dvi and ps output
  \ifthenelse{\boolean{pdfoutput}}%
	     {\typeout{using pdf}
	       \usepackage[pdftex]{graphicx,xcolor}
	       \newcommand{\graphext}{pdf}
	       \newcommand{\graphextex}{pdf_t}
	       \usepackage{epsfig}
	       \usepackage{tikz}
	     }
	     {
	       \typeout{using eps}
	       \usepackage[dvips]{graphicx,xcolor}
	       \newcommand{\graphext}{eps}
	       \newcommand{\graphextex}{eps_t}
	       \usepackage{epsfig}
	       \usepackage{tikz}
	     }
}
\newcommand{\hoz}{\sobhz1(\W)}

\newcommand{\T}[1]{\cT^{#1}}

\renewcommand{\ie}{i.e.,\xspace}

\providecommand{\figwidth}{\textwidth} \newcommand{\figscale}{1}

\renewcommand{\fes}{\ensuremath{\mathbb V}}

\newcommand{\MA}{Monge-Amp\`ere }
\newcommand{\cof}[1]{\operatorname{cof} #1}
\renewcommand{\det}[1]{\operatorname{det} #1}
\newcommand{\tensorp}[2]{\ensuremath{\left\langle#1\otimes#2\right\rangle}}
\newcommand{\frob}[2]{\ensuremath{{#1}{:}{#2}}}
\renewcommand{\div}[1]{\operatorname{div}\left(#1\right)}
\renewcommand{\H}{\ensuremath{\vec{H}}}

\newcommand{\symm}{\ensuremath{\operatorname{Sym}(\reals^{d\times d})}}

\newcommand{\A}{\ensuremath{\vec{A}}}

\numberwithin{equation}{section}
\setboolean{showtodo}{false}
\setboolean{shownotes}{false}
\setboolean{showchanges}{false}
\setboolean{usemathrsfs}{false}%true
\pdfformat{true}%false

\newcommand{\fenics}{\Program{FEniCS}}
\newcommand{\dolfin}{\Program{DOLFIN}}

\renewcommand{\H}{\ensuremath{\geovec{H}}}

\newcommand{\MAD}{Monge--Amp\`ere--Dirichlet }

\author{Tristan Pryer} 
\address{ Tristan Pryer\newline
  Department of Mathematics, Statistics and Actuarial Science\newline 
  University of Kent\newline
  Canterbury\newline 
  CT2 7NF, United Kingdom} 
\curraddr{}
\email{\linkedemail{T.Pryer@kent.ac.uk}}

\title[Applications of NVFEMs to MAD equations]{Applications of nonvariational finite element methods to Monge--Amp\`ere type equations} 
\date{\today}
\begin{document}
\maketitle

\begin{abstract}The goal of this work is to illustrate the application of
  the nonvariational finite element method to a specific \MA
  type nonlinear partial differential equation. The equation we
  consider is that of prescribed Gauss curvature although the method
  can be generalised to any \MA operator.
\end{abstract}

\section{Introduction and problem setting}
\label{sec:intro}

The \emph{nonvariational finite element method} (NVFEM) introduced in
\cite{LakkisPryer:2011} is a numerical method aimed at problems of the
form
\begin{equation}
  \label{eq:nonvar-op}
  \frob{\A(\vec x)}{\Hess u(\vec x)} = f(\vec x)
\end{equation}
where for each $\vec x \in \W$ the matrix $\A(\vec x)\in\symm$, the space of
bounded symmetric positive definite matrices and where $\Hess u$
denotes the Hessian of the function $u$. The operation
$\frob{\geovec B}{\geovec C} = \trace\qp{\Transpose{\geovec B}\geovec C}$
is the Frobenius inner product between two $d\times d$
matrices. Classical finite element methods are applicable to this
problem if we assume the coefficient matrix $\A$ is differentiable. In
this case we may rewrite (\ref{eq:nonvar-op}) in \emph{variational} or
\emph{divergence} form via the introduction of an advection term since
\begin{equation}
  \label{eq:var-op}
  f(\vec x)
  =
  \frob{\A(\vec x)}{\Hess u(\vec x)} 
  =
  \div{\qp{\A(\vec x)\nabla u(\vec x)}} - \D{\A(\vec x)}{\nabla u(\vec x)}
\end{equation}
where 
\begin{equation}
  \D\A(\geovec x) 
  =
  \qp{\sum_{i=1}^d \partial_i a_{i,1}(\geovec x) , \dots , \sum_{i=1}^d \partial_i a_{i,d}(\geovec x)}.
\end{equation}
Note that we are using the convention that $\nabla \phi =
\Transpose{\qp{\partial_1\phi, \dots,\partial_d\phi}}$ is the column
vector formed of first order partial derivatives of a $d$--multivariate
function $\phi$.

The introduction of the advection term may result in the variational
problem becoming advection dominated. This is undesirable in the
finite element context and stabilisation terms become necessary to
derive a viable numerical method
\cite[c.f.]{ErnGuermond:2004}. Interestingly if
$\Norm{\D\A}_{\leb{\infty}(\W)} \gg \Norm{\A}_{\leb{\infty(\W)}}$
applying the NVFEM to (\ref{eq:nonvar-op}) does not result in an
unstable scheme, whereas applying a standard FEM to (\ref{eq:var-op})
does. This is numerically demonstrated in \cite[\S
  4.2]{LakkisPryer:2011}. It may even be the case that $\A$ is not
differentiable, in which case the standard FEM cannot be applied. 

The fully nonlinear problem
\begin{equation}
  \label{eq:nonlin-op}
  \cF(\Hess u) = 0
\end{equation}
is related to the nonvariational problem (\ref{eq:nonvar-op}) by the
fundamental theorem of calculus. If $\cF$ is sufficiently regular and
$u$ solves (\ref{eq:nonlin-op}) then $u$ also solves
\begin{equation}
  \frob{\left[ \int_0^1 \cF'(t \Hess u)  \d t  \right]}{\Hess u} + \cF(0) 
  =
  0,
\end{equation}
where $\cF'$ is the Fr\'echet derivative of $\cF$. It is also the case
that a Newton linearisation applied to (\ref{eq:nonlin-op}) results in
a sequence of linear nonvariational PDEs:
\begin{equation}
  \frob{\cF'(\Hess u^n)}{\Hess \qp{u^{n+1}-u^n}} = -\cF(\Hess u^n).
\end{equation}

%% \begin{comment}
%% In this work we are interested in a nonlinear fourth order PDE
%% \begin{equation}
%%   \label{eq:4thorder}
%%   \frob{\cof\qp{\Hess u}}{\Hess\qp{\det\qp{\Hess u}^{\frac{-\qp{d+1}}{d+2}}}} = 0,
%% \end{equation}
%% where $\cof\qp{\vec B}$ denotes the matrix of cofactors of $\vec
%% B$. This arises in the following setting. Let $u$ be a strictly convex
%% function then it follows that $\det{\Hess u} >0$ and we may introduce
%% a uniformly convex hypersurface $\cM$ with Gaussian curvature
%% $\kappa$. The affine area functional typically denoted in the literature as
%% \begin{equation}
%%   \cA[u] := \int_\W \qp{\det \Hess u}^{\frac{1}{d+2}} = \int_{\cM_\W} \kappa^{\frac{1}{d+2}},
%% \end{equation}
%% where
%% \begin{equation}
%%   \cM_\W := \{ \qp{\vec x, u(\vec x)} \in \cM | \vec x\in\W \}
%% \end{equation}
%% is the graph of the function $u$ in $\cM$.

%% The PDE (\ref{eq:4thorder}) can be written as a coupled system of two
%% second order PDEs, in this case we seek a pair $u, v$ such that
%% \begin{gather}
%%   \label{eq:cof-diff}
%%   \frob{\cof{\Hess u}}{\Hess v} = 0
%%   \\
%%   \label{eq:det-diff}
%%   \det\qp{\Hess u}^{\frac{-\qp{d+1}}{d+2}} - v = 0.
%% \end{gather}
%% The first, (\ref{eq:cof-diff}), being a linear PDE in $v$, the second,
%% (\ref{eq:det-diff}), a fully nonlinear PDE in $u$. This pair of PDEs
%% closely resemble Monge Ampere type equations, indeed
%% (\ref{eq:det-diff}) is of Monge Ampere type whereas
%% (\ref{eq:cof-diff}) is the linearisation of a Monge Ampere operator.
%% \end{comment}

The \MA operators are an extremely interesting class of fully
nonlinear PDE. These arise from differential geometry and optimal
transport problems; they take the form
\begin{equation}
  \label{eq:MA-operators}
  \cF(\Hess u):=\det\qp{\Hess u} - f\qp{\nabla u, u, \geovec x} = 0.
\end{equation}
For example the \MAD (MAD) problem is the case when $f = f(\geovec x)$
and (\ref{eq:MA-operators}) is coupled with a Dirichlet type boundary
condition ($u = g$ on $\partial \W$). This particular equation is a
prototypical example of a fully nonlinear PDE. %Numerical studies on
%this problem include (but are not limited to)
%\cite{DeanGlowinski:2006, OlikerPrussner:1988, LoeperRapetti:2005, publication7}.

There are a variety of numerical methods available 
for the more
general \MA class of fully nonlinear PDE (\ref{eq:MA-operators}). In
\cite{Oberman:2008} the author proposes a wide stencil finite
difference scheme. In \cite{Bohmer:2008} a $\cont{1}$ finite element
scheme based on the Argyris element is used. In a series of papers
Feng and Neilan \cite{FengNeilan:2009,FengNeilan:2009a} construct
numerical approximations of solutions to sequences of quasilinear
biharmonic equations. This is very reminiscent of the vanishing
viscosity method first studied for use in fully nonlinear first order
PDEs. The method is aptly named the vanishing moment method. More
recently in \cite{BrennerGudiNeilanSung:2011} a consistent
penalisation method has been introduced for these problems. Finally,
Awanou \cite{Awanou:2011} uses a \emph{Laplacian relaxation} technique
to study these equations.

For the \MA type equation (\ref{eq:MA-operators}) to be well
posed we require $\W\subset \reals^d$ to be a convex domain and
$f>0$. The \MA operator is elliptic over the cone of strictly
convex functions in $\W$ and under the constraints above will admit a
unique convex viscosity solution \cite{CaffarelliCabre:1995}.

In this work we will study the equation of prescribed Gauss
curvature. This arises from the problem of finding a function $u$ such
that the graph of $u$ has a specified Gaussian curvature $K$. In this
case we have that
\begin{equation}
  K = \frac{\det{\Hess u}}{\qp{1+\norm{\nabla u}^2}^{(d+2)/2}}
\end{equation}
and hence 
\begin{equation}
  \label{eq:prescribed-gauss-curvature}
  \cF(\Hess u, \nabla u, \geovec x) 
  :=
  \det{\Hess u} - K \qp{1+\norm{\nabla u}^2}^{(d+2)/2}.
\end{equation}
Note that $K = K(u, \geovec x)$.

The linearisation of this problem can be calculated in a direction $v$
as
\begin{equation}
  \label{eq:linearisation}
  \begin{split}
    \frob{\cF'(\Hess u, \nabla u, \geovec x)}{\Hess v}
    &=
    \lim_{\epsilon \to 0} 
    \frac{1}{\epsilon}
    \qp{
      \cF(\Hess u + \epsilon\Hess v, \nabla u + \epsilon\nabla v, \geovec x) 
      - 
      \cF(\Hess u, \nabla u, \geovec x) 
    }
    \\
    &=
    \frob{\cof{\Hess u}}{\Hess v} 
    +
    \qp{d+2}K\qp{\qp{1+\norm{\nabla u}^2}^{d/2}\Transpose{\qp{\nabla u}}\nabla v}
  \end{split}
\end{equation}
and thus the linearisation is elliptic if $\cof{\Hess u}$ is an
elliptic operator. This holds for convex $u$.

\section{Discretisation}

The process of discretisation can be sought in two ways. We may look
at the continuous problem and discretise this directly, resulting in a
system of nonlinear equations, or we may first linearise the problem
and discretise from there. Discretising the nonlinear problem directly
is certainly possible but is more technical, as discussed in
\cite{BrennerGudiNeilanSung:2011}. For brevity we will perform a
Newton linearisation on (\ref{eq:MA-operators}) and discretise the sequence of
linear nonvariational PDEs in a similar light to \cite{LakkisPryer:2011b}.

Let $\T{}$ be a conforming, shape regular triangulation of $\W$,
namely, $\T{}$ is a finite family of sets such that
\begin{enumerate}
\item $K\in\T{}$ implies $K$ is an open simplex (segment for $d=1$,
  triangle for $d=2$, tetrahedron for $d=3$),
\item for any $K,J\in\T{}$ we have that $\closure K\meet\closure J$ is
  a full sub-simplex (i.e., it is either $\emptyset$, a vertex, an
  edge, a face, or the whole of $\closure K$ and $\closure J$) of both
  $\closure K$ and $\closure J$ and
\item $\union{K\in\T{}}\closure K=\closure\W$.
\end{enumerate}
We also define $\cE$ to be the skeleton of the triangulation, that is
the set of sub-simplices of $\T{}$ contained in $\W$ but not $\partial
\W$. For $d=2$, for example, $\cE$ would consist of the set of edges of
$\T{}$ not on the boundary.

We use the convention where $\funk h\W\reals$ denotes the
\emph{meshsize function} of $\T{}$, i.e.,
\begin{equation}
  h(\vec{x}):=\max_{\closure K\ni \vec x}h_K.
\end{equation}

\begin{Defn}[FE spaces]
   Let $\poly p(\T{})$ denote the space of piecewise polynomials of
   degree $k$ over the triangulation $\T{}$ of $\W$. We introduce the
   \emph{finite element spaces}
   \begin{gather}
     \label{eqn:def:finite-element-space}
     \fes{} = \poly p(\T{}) \cap \cont{0}(\W) \cap \hoz \AND
     \fesW{} = \poly p(\T{}) \cap \cont{0}(\W)
   \end{gather}
   to be the usual space of continuous piecewise polynomial functions and
   \begin{equation}
     \rS := \fes{} \times \fesW^{d^2}.
   \end{equation}
 \end{Defn}
 
\begin{Rem}[generalised Hessian]
  \label{rem:generalised-hessian}
  Given a function $v\in\sobh2(\W)$, let $\geovec
  n:\partial\W\to\reals^d$ be the outward pointing normal of $\W$ then the
  Hessian of $v$, $\Hess v$, satisfies the following identity:
  \begin{equation}
    \label{eq:generalised-hessian}
    \duality{\Hess v}{\phi}
    = 
    -
    \int_\W{\nabla v}\otimes{\nabla \phi}
    +
    \int_{\partial\W}{\nabla v}\otimes{\geovec n \ \phi}
    \Foreach \phi\in\sobh1(\W)
  \end{equation}
  where $\duality{\cdot}{\cdot}$ denotes an appropriate duality
  pairing. It follows that we may weaken the regularity assumptions
  to $v\in\sobh1(\W)\cap\sobh1(\partial\W)$.
\end{Rem}

\begin{Defn}[finite element Hessian]
  \label{defn:feh}
  From Remark \ref{rem:generalised-hessian} and in view of Reisz
  representation theorem we may define the \emph{finite element
    Hessian} such that
  \begin{equation}
    \int_\W {\H[V]}{\Phi} = \duality{\Hess V}{\Phi} \Foreach \Phi\in\fesW{}.
  \end{equation}
\end{Defn}

\begin{Pro}[symmetry of the finite element Hessian]
  \label{prop:symm}
  The finite element Hessian is symmetric, that is for each
  $V\in\fes{}$
  \begin{equation}
    \int_\W \H[V] \Phi = \int_\W \Transpose{\qp{\H[V]}} \Phi
    \Foreach \Phi\in\fesW{}.
  \end{equation}
\end{Pro}

In view of the constraints to the continuous problem
(\ref{eq:prescribed-gauss-curvature}) to admit a unique solution it is
also necessary to construct a discrete notion of convexity. This has
been developed in \cite{Aguilera:2008} and is naturally passed down
from the concept of distributional convexity.

%% \begin{theorem}[consistency of the FE hessian \cite{Pryer:2010}]
%%   Let the conditions of Theorem \ref{the:inf-sup} hold. In addition
%%   assume that $u\in\sobh{j}(\W)$ for some $j\in [0,p+1]$, then we have
%%   that there exists a $C>0$ such that
%%   \begin{equation}
%%     \Norm{\Hess u - \H[U]}_{\sobh{-1}(\W)}
%%     \leq 
%%     C h^{j+1} \norm{u}_j.
%%   \end{equation}
%%   This immediately infers that 
%%   \begin{equation}
%%     \Norm{\Hess u^n - \H[U^n]}_{\sobh{-1}(\W)}
%%     \leq 
%%     C h^{j+1} \norm{u^n}_j.
%%   \end{equation}
%% \end{theorem}

\begin{Defn}[finite element convexity \cite{Aguilera:2008}]
  \label{defn:fe-convex}
  A function, $v\in\sobh1(\W)\cap\sobh1(\partial\W)$, is said to be
  \emph{finite element convex} if
  \begin{equation}
    \label{eq:fe-convex}
    \int_\W{\H[v]}{\Phi} \text{ is positive semidefinite } \Foreach \Phi\in\fesW{}
  \end{equation}
  where $\Phi \geq 0$ on $\W$. It is strictly finite element convex if
  (\ref{eq:fe-convex}) is positive definite.
\end{Defn}

Definition \ref{defn:feh} allows us to construct what is essentially a
2--0 mixed method where the Hessian of the solution to
(\ref{eq:MA-operators}) is treated as an auxiliary variable in the
formulation (as opposed to the 1--1 mixed methods commonly found in
the literature by decoupling a second order PDE into a system of first
order PDEs \cite[c.f.]{BrezziFortin:1991}). Given the linearisation
(\ref{eq:linearisation}) we formulate the problem in the discrete
setting as follows: Given an initial guess $\qp{U^0, \H[U^0]} \in \rS$
that is strictly finite element convex (\ref{eq:fe-convex}), for
$n\in\naturals$ find $\qp{U^{n}, \H[U^{n}]} \in \rS$ such that
\begin{gather}
  \label{eq:hessian-representation}
  \ltwop{\H[U^{n}]}{\Phi} 
  +
  \tensorp{\nabla U^{n}}{\nabla \Phi}
  -
  \sum_{e\in\partial\W} 
  \tensorp{\nabla U^{n}}{\geovec n \ \Phi}_e
  =
  0
  \\
  \label{eq:nonlinear-representation}
  \ltwop{\frob{\cF'\qp{\H[U^{n-1}], \nabla U^{n-1}}}{\H[U^{n}-U^{n-1}]} + \cF\qp{\H[U^{n-1}],\nabla U^{n-1}}}{\Phi}
  = 
  0
  \Foreach \Phi\in\fes{}.
\end{gather}
Where for the problem of prescribed Gaussian curvature (\ref{eq:nonlinear-representation}) is
%% \begin{equation}
%%   \begin{split}
%%     \ltwop{\frob{\cF'\qp{\H[U^{n-1}], \nabla U^{n-1}}}{\H[U^{n}-U^{n-1}]} + \cF\qp{\H[U^{n-1}]},\nabla U^{n-1}}{\Phi}
%%     \\
%%     =
%%     \int_\W \frob{\cof{\H[U^{n-1}]}}{\H[U^{n}-U^{n-1}]} \Phi&
%%     \\
%%     +
%%     \int_\W 2dK {\qp{1+\norm{\nabla U^{n-1}}^2}^{d/2}\nabla U^{n-1}}\nabla U^n\Phi&
%%     \\
%%     + 
%%     \int_\W \qp{\det{\H[U^{n-1}]} - K \qp{1+\norm{\nabla U^{n-1}}^2}^{(d+2)/2}}\Phi&.
%%   \end{split}
%% \end{equation}
\begin{eqnarray}
  \label{eqn:linearisation}
    &&\ltwop{\frob{\cF'\qp{\H[U^{n-1}], \nabla U^{n-1}}}{\H[U^{n}-U^{n-1}]} + \cF\qp{\H[U^{n-1}],\nabla U^{n-1}}}{\Phi}
    \nonumber \\
    && \qquad \qquad \qquad \qquad 
    = 
    \int_\W \frob{\cof{\H[U^{n-1}]}}{\H[U^{n}-U^{n-1}]} \Phi \\
    && \qquad \qquad \qquad \qquad \qquad
    +
    \int_\W 2dK {\qp{1+\norm{\nabla U^{n-1}}^2}^{d/2}\Transpose{\qp{\nabla U^{n-1}}}}\nabla \qp{U^{n}-U^{n-1}} \Phi \nonumber
    \\
    && \qquad \qquad \qquad \qquad  \qquad 
     + 
    \int_\W \qp{\det{\H[U^{n-1}]} - K \qp{1+\norm{\nabla U^{n-1}}^2}^{(d+2)/2}}\Phi      \nonumber.
\end{eqnarray}

Due to the symmetry property given in Proposition \ref{prop:symm} we
may simplify the problem somewhat to seeking only the upper (or lower)
triangular parts of the finite element Hessian. This reduces $\rS =
\fes{} \times \fesW^{(d^2 + d)/2}$.

\begin{The}[solvability of the discrete system \cite{Pryer:2010}]
  \label{the:inf-sup}
  Let $U \in \fes{}$ be the nonvariational finite element
  approximation to $u$, the solution of
  \begin{equation}
    \frob{\A}{\Hess u} = f,
  \end{equation}
  where $\A$ is an elliptic operator. Then we have a discrete inf--sup
  condition, that is the linear system is always invertible. Hence,
  assuming the linearisation maintains ellipticity, the discrete
  problem
  (\ref{eq:hessian-representation})--~(\ref{eq:nonlinear-representation})
  is well posed.
\end{The}

\section{Numerical experiments}

In this section we detail numerical experiments on the formulation
(\ref{eq:hessian-representation})--~(\ref{eq:nonlinear-representation}). 

We will consider the case $d=2$ and when $K > 0$ is some prescribed
curvature. In each of the experiments we choose $p=2$, \ie $\fes{}$
consists of piecewise quadratic functions. The domain $\W$ is taken as
a square whose size differs on each of the experiments and the
triangulation $\T{}$ is unstructured. All of the numerical experiments
have been conducted using the \dolfin environment of the finite
element package \fenics \cite{LoggWells:2010}.

In Figures
\ref{fig:convergence-rates}--\ref{fig:convergence-rates-exp} we
construct classical solutions to (\ref{eq:prescribed-gauss-curvature})
in order to look at the numerical convergence of the method.  In
Figure \ref{fig:various-Ks} we consider $K$ as a constant over the
domain $[-.57,.57]^2$. These results can then be compared with the two
other numerical studies found in the literature on prescribed Gauss
curvature \cite{FengNeilan:2009, Awanou:2011}. In these experiments
the authors note that the problem
(\ref{eq:prescribed-gauss-curvature}) is well posed only for $K \leq
K^{\max{}}$ and estimate the value of $K^{\max{}}$ by asserting when the
numerical algorithm proposed breaks down.

The initial guess to any Newton iteration is paramount due to the well
known \emph{overshoot} property. In the case of \MA type
linearisations it's especially important since (discrete) convexity
must be maintained during the iterative procedure for the problem to
remain well posed. In each of the tests below we initialise the
algorithm by approximating the solution of the MAD problem over the
initial mesh as detailed in \cite{LakkisPryer:2011b}.
\begin{equation}
  \begin{split}
    \det{\Hess u} &= K \quad\text{ in } \W
    \\
    u &= g \quad\text{ on } \partial\W.
  \end{split}
\end{equation}

The iterative procedure given by the discrete problem
(\ref{eq:hessian-representation})--~(\ref{eq:nonlinear-representation})
is terminated when two concurrent iterates satisfy $\Norm{U^n -
  U^{n-1}}_{\leb{\infty}(\W)} \leq 10^{-10}$.

%%%%%%%%%%%%%%%%%%%%%%%%%%%%%%%%%%%%%%%%%%%%%%%%%%%%%%%%%%%%%%%%%%%%%%%%
%%%%%%%%%%%%%%%%%%%%%%%%%%%%%%%%%%%%%%%%%%%%%%%%%%%%%%%%%%%%%%%%%%%%%%%%
%%%%%%%%%%%%%%%%%%%%%%%%%%%%%%%%%%%%%%%%%%%%%%%%%%%%%%%%%%%%%%%%%%%%%%%%
%  FIGURES!!
%%%%%%%%%%%%%%%%%%%%%%%%%%%%%%%%%%%%%%%%%%%%%%%%%%%%%%%%%%%%%%%%%%%%%%%%
%%%%%%%%%%%%%%%%%%%%%%%%%%%%%%%%%%%%%%%%%%%%%%%%%%%%%%%%%%%%%%%%%%%%%%%%
%%%%%%%%%%%%%%%%%%%%%%%%%%%%%%%%%%%%%%%%%%%%%%%%%%%%%%%%%%%%%%%%%%%%%%%%
\renewcommand{\figscale}{0.29}

\begin{figure}[!htbp]
  \caption{\label{fig:convergence-rates} In this experiment we fix
    choose a convex solution $u$ which classically solves the equation
    of prescribed Gauss curvature
    (\ref{eq:prescribed-gauss-curvature}) over the square
    $[-.5,.5]^2$. That is, we fix $u = \norm{\geovec x}^4$ and
    calculate $K = K(\geovec x, u)$. We solve the discrete problem
    over a sequence of concurrently refined meshes and ascertain the
    errors and convergence rates for the problem in $\leb{2}(\W)$,
    $\hoz$ and a discrete $\sobh{2}(\W)$ seminorm. Notice that
    $\Norm{u-U^N} \approx \Oh(h^3)$, $\norm{u - U^N}_1 \approx
    \Oh(h^2)$ and $\Norm{\Hess u - \H[U^N]}\approx \Oh(h^{1.5})$.}
  \begin{center}
    \subfigure[Errors and convergence rates for the
      problem, $p=2$.] {
      \includegraphics[scale=\figscale,width=0.47\figwidth]{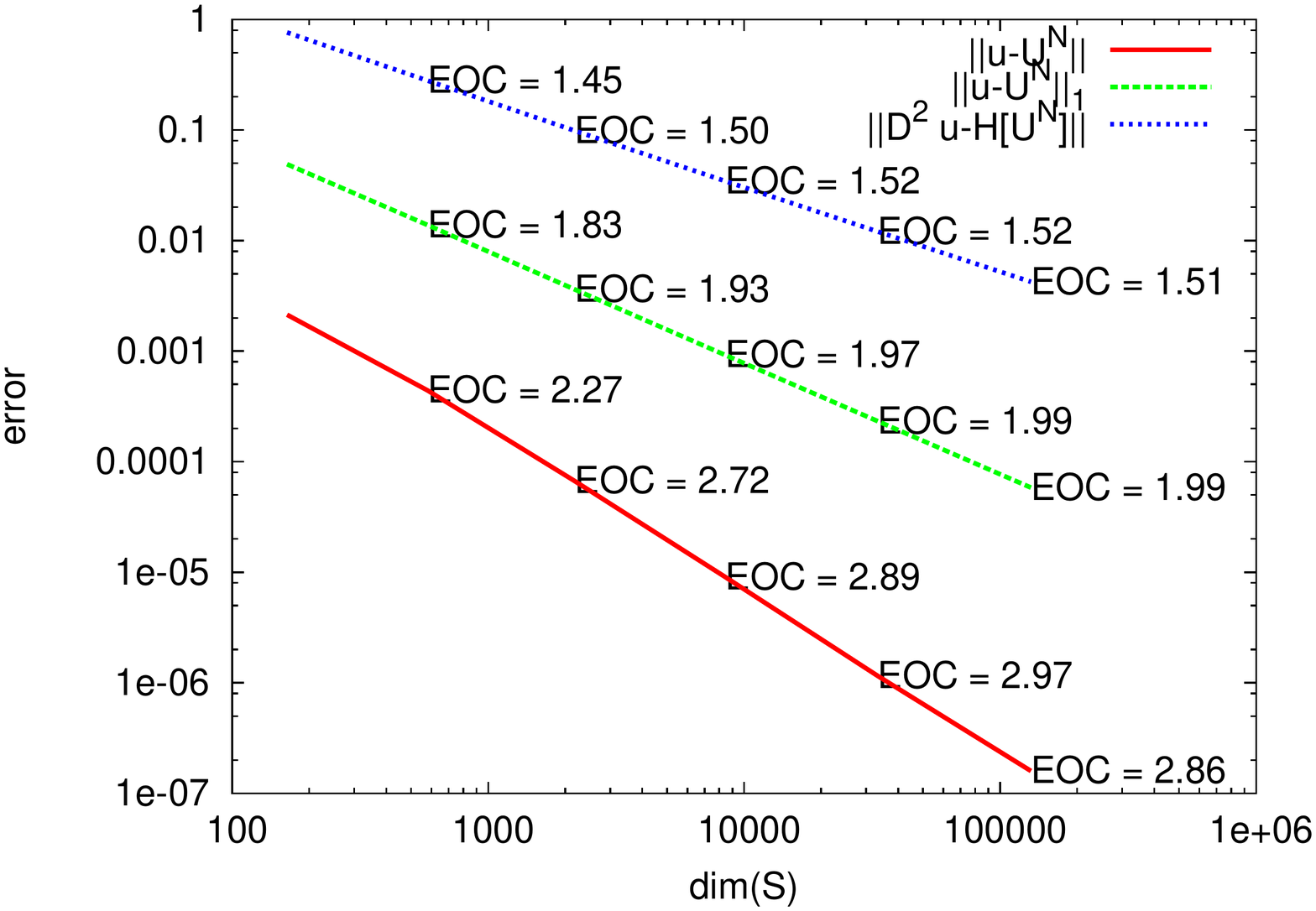}
    }  \hfill
    \subfigure[Solution plot] {
      \includegraphics[scale=\figscale,width=0.47\figwidth]{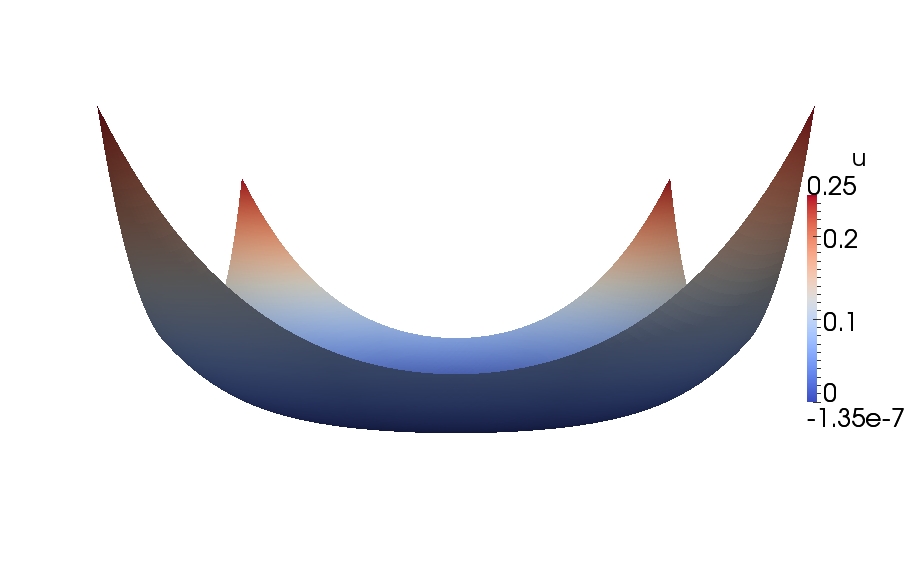}
    }
  \end{center}
\end{figure}

\begin{figure}[!htbp]
  \caption[]{\label{fig:convergence-rates-exp} In this experiment we fix
    choose a convex solution $u$ which classically solves the equation
    of prescribed Gauss curvature
    (\ref{eq:prescribed-gauss-curvature}) over the square
    $[-.5,.5]^2$. That is, we fix $u = \exp\qp{\norm{\geovec x}^2/2}$ and
    calculate $K = K(\geovec x, u)$. We solve the discrete problem
    over a sequence of concurrently refined meshes and ascertain the
    errors and convergence rates for the problem in $\leb{2}(\W)$,
    $\hoz$ and a discrete $\sobh{2}(\W)$ seminorm. Notice that
    $\Norm{u-U^N} \approx \Oh(h^3)$, $\norm{u - U^N}_1 \approx
    \Oh(h^2)$ and $\Norm{\Hess u - \H[U^N]}\approx \Oh(h^{1.5})$.}

  \subfigure[Errors and convergence rates for the
    problem, $p=2$.] {
    \includegraphics[scale=\figscale,width=0.47\figwidth]{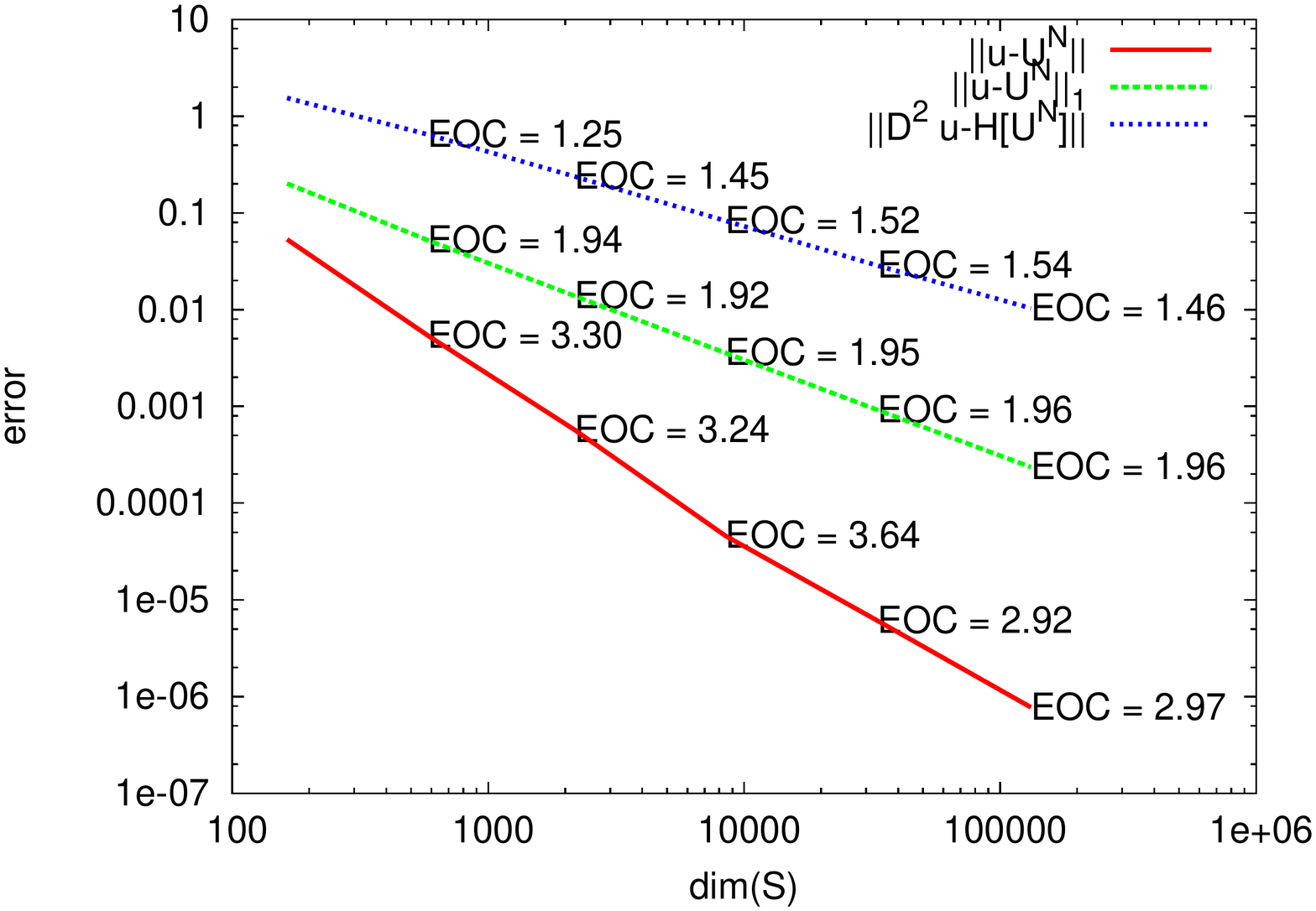}
  }  \hfill
  \subfigure[Solution plot] {
    \includegraphics[scale=\figscale,width=0.47\figwidth]{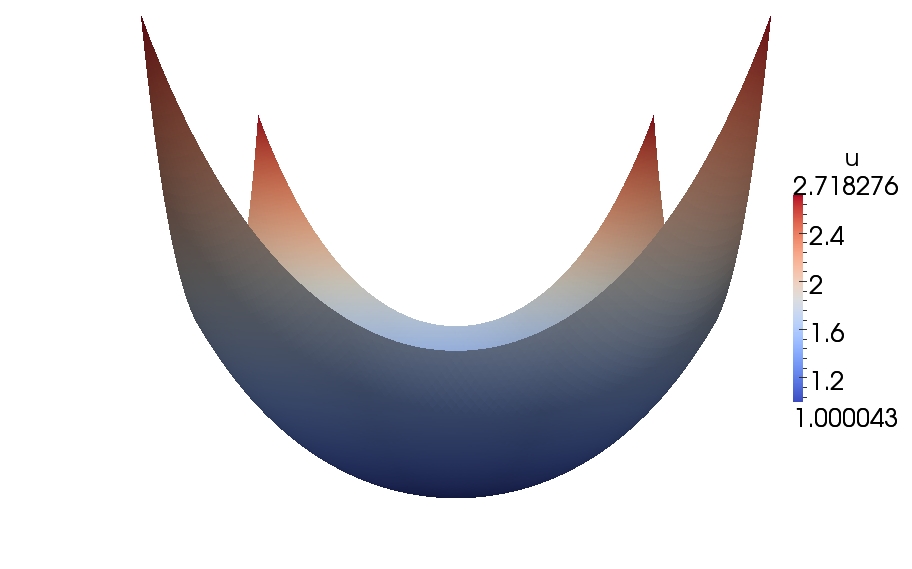}
  }
\end{figure}

\begin{figure}[!htbp]
  \caption[]{\label{fig:various-Ks} In this experiment we fix $h
    \approx 0.009$ and implement the discrete problem over an
    unstructured mesh of the square $[-0.57,0.57]^2$ consider the case
    $K$ is constant. We choose various values of $K>0$ and display a
    contour plot together with a side view of the discrete
    solution. Note that the numerical algorithm fails to converge for
    $K=2$.}  \subfigure[Contour plot for $K=0.01$] {
    \includegraphics[scale=\figscale,width=0.45\figwidth]{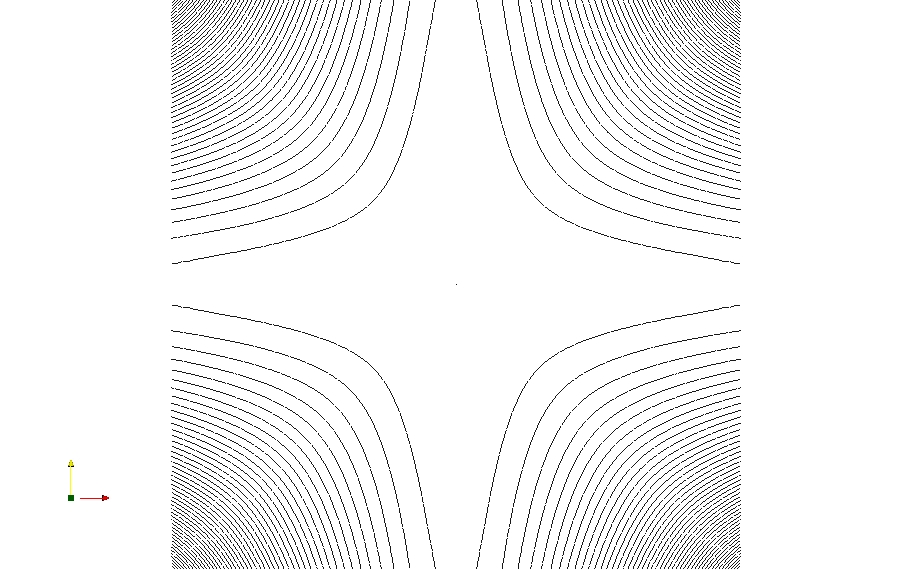}
  } \hfill \subfigure[Solution plot for $K=0.01$] {
    \includegraphics[scale=\figscale,width=0.45\figwidth]{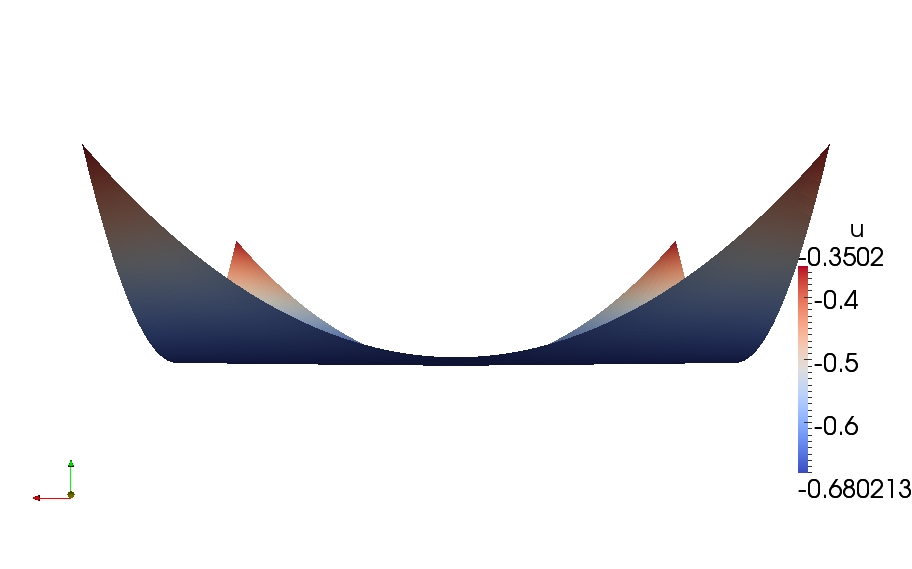}
  } \subfigure[Contour plot for $K=0.1$] {
    \includegraphics[scale=\figscale,width=0.45\figwidth]{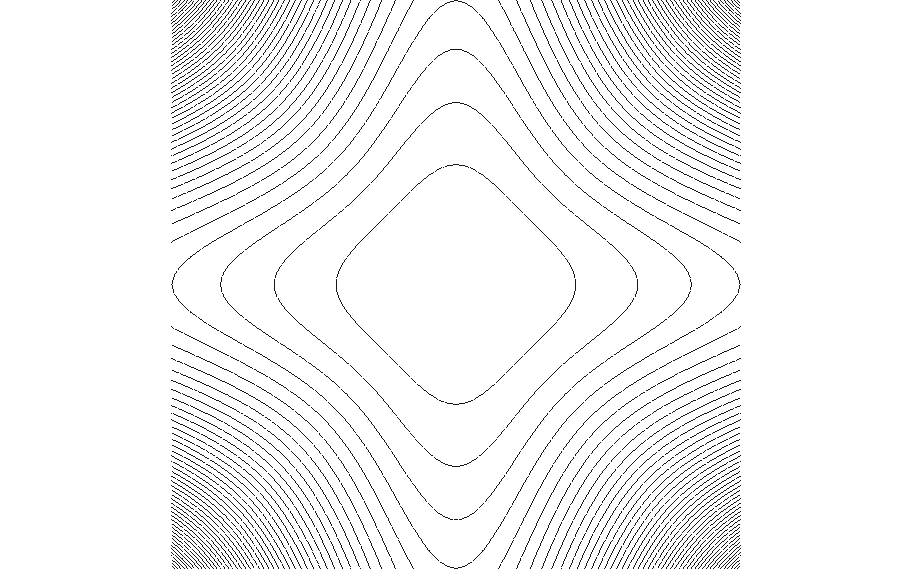}
  } \hfill \subfigure[Solution plot for $K=0.1$] {
    \includegraphics[scale=\figscale,width=0.45\figwidth]{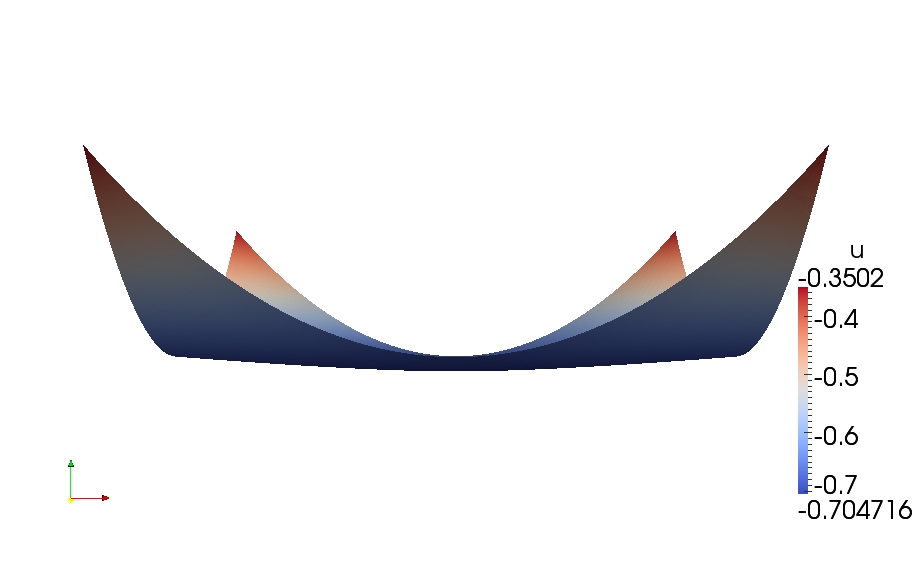}
  } \subfigure[Contour plot for $K=0.5$] {
    \includegraphics[scale=\figscale,width=0.45\figwidth]{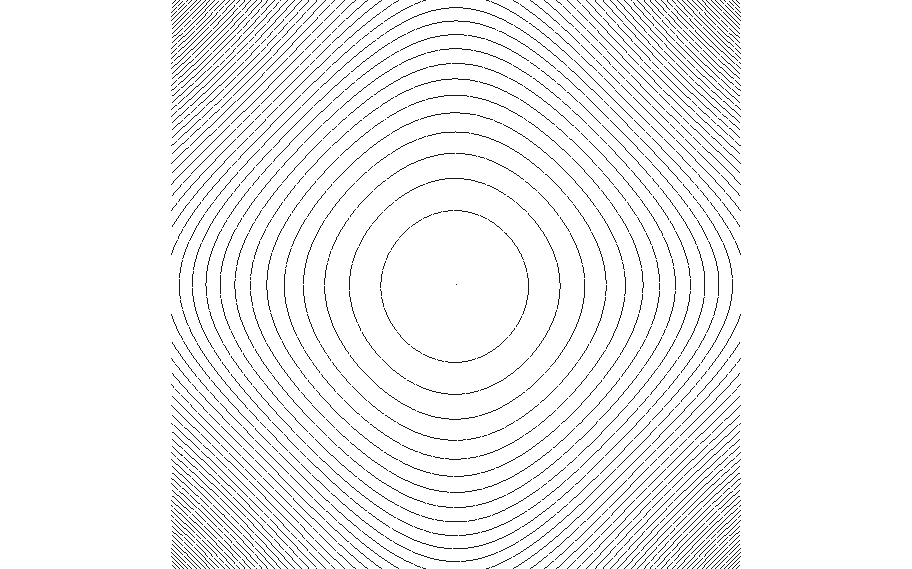}
  } \hfill \subfigure[Solution plot for $K=0.5$] {
    \includegraphics[scale=\figscale,width=0.45\figwidth]{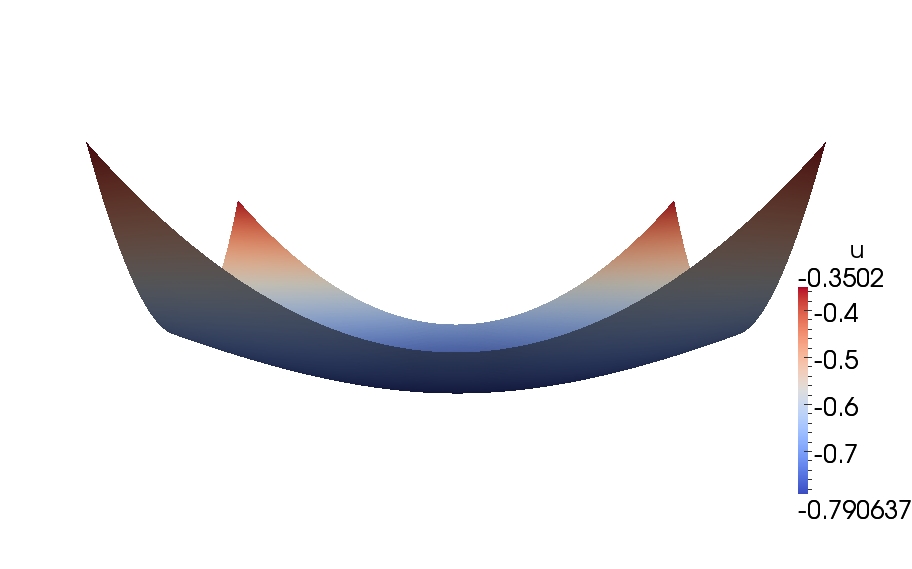}
  } \subfigure[Contour plot for $K=1.0$] {
    \includegraphics[scale=\figscale,width=0.45\figwidth]{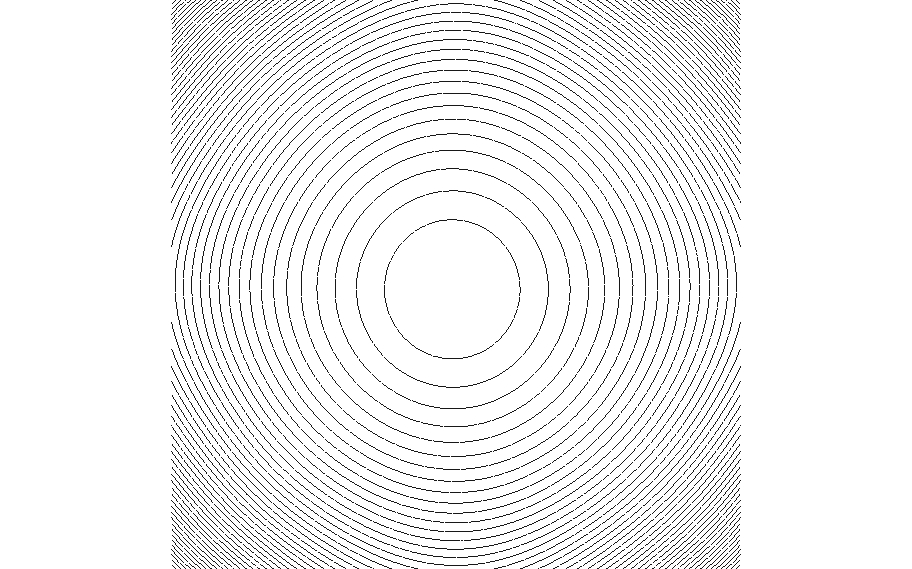}
  } \hfill \subfigure[Solution plot for $K=1.0$] {
    \includegraphics[scale=\figscale,width=0.45\figwidth]{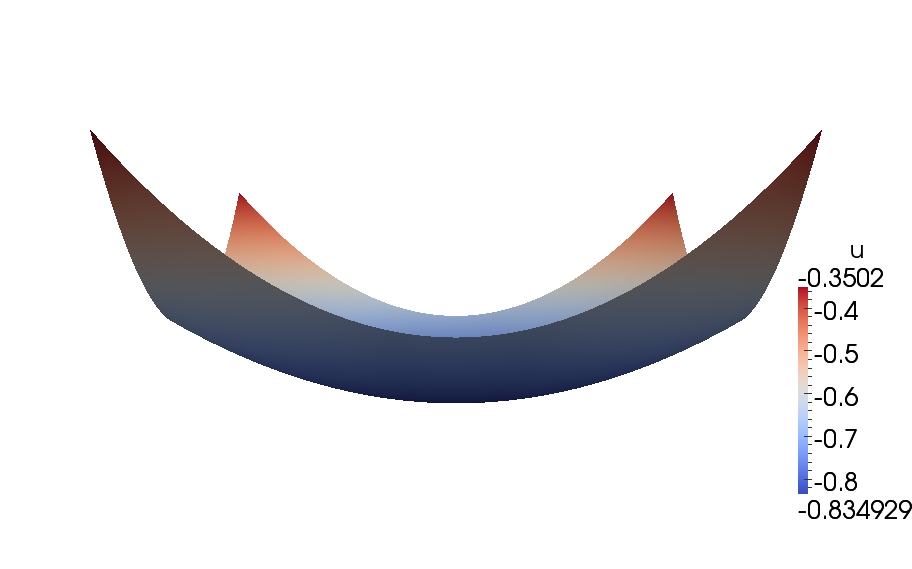}
  } \subfigure[Contour plot for $K=1.5$] {
    \includegraphics[scale=\figscale,width=0.45\figwidth]{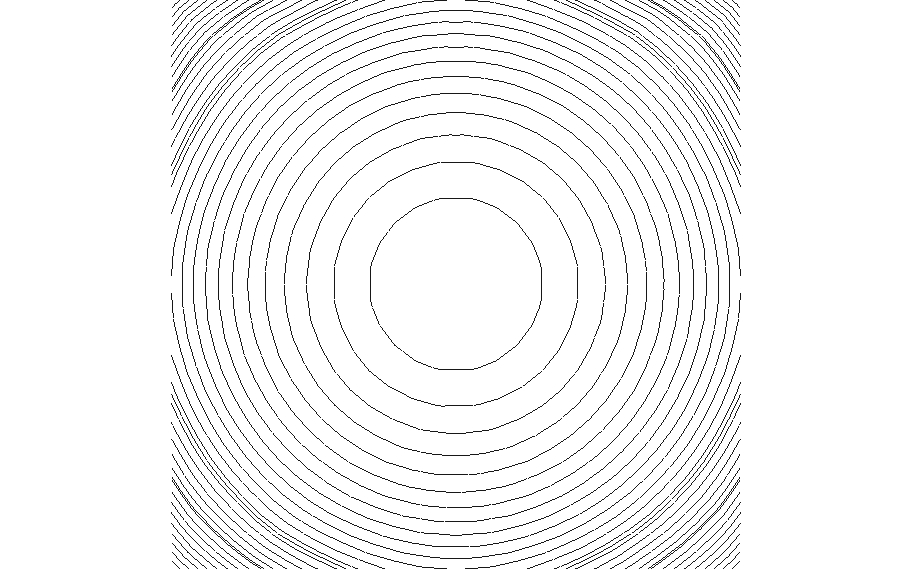}
  } \hfill \subfigure[Solution plot for $K=1.5$] {
    \includegraphics[scale=\figscale,width=0.45\figwidth]{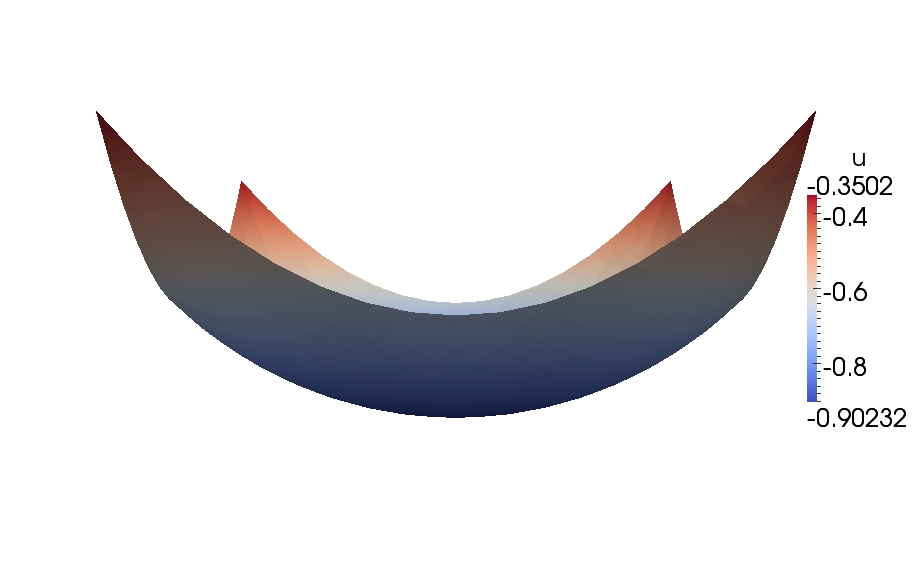}
  }
\end{figure}

\bibliographystyle{alpha}
\bibliography{tristansbib}

\end{document}